\documentclass[11pt]{elsarticle}\pdfoutput=1
\usepackage[T1]{fontenc}\usepackage[latin1]{inputenc}
\usepackage{times}
\usepackage{parskip}\usepackage[normalem]{ulem}
\usepackage[colorlinks,linkcolor=black,citecolor=black,urlcolor=black,%
            bookmarks=false]{hyperref}
\usepackage{amsmath,amsfonts,amssymb,amsopn}
\usepackage[mathscr]{euscript}
\usepackage{bbold}
\makeatletter\let\@EA\expandafter
\let\o@Lim\lim\def\lim{\o@Lim\limits}
\let\o@Sum\sum\def\sum{\o@Sum\limits}
\let\v@rphi\phi\let\phi\varphi\let\varphi\v@rphi
\def\@endproof{\QED\endtrivlist}
\def\@ddtok#1\@{\@temptokena\@EA{\the\@temptokena#1}}
\def\@tokX#1=#2\@#3{\@EA\@ddtok\csname#3#1\endcsname\@%
  \ifx\\#2\\\@ddtok{{#1}}\@\else\@tokY#2\fi}
\def\@tokY#1={\@ddtok{{#1}}\@}
\def\MathSyms#1{\@temptokena{}\@dblarg{\M@thSym{#1}}}
\def\M@thSym#1[#2]#3{\def\m@tSX##1##2##3{\def##2{##1{##3}}}\def\m@tSY{%
    \@temptokena{\m@tSX}\@EA\@ddtok\csname math#1\endcsname\@\@tokX}%
  \@for\@gtempa:=#3\do{\@EA\m@tSY\@gtempa=\@{#2}\the\@temptokena}}
\def\MathOps#1{%
  \def\m@tO{\@temptokena{\DeclareMathOperator}\@ifstar{\@ddtok*\@\@tokX}\@tokX}%
  \@for\@gtempa:=#1\do{\@EA\m@tO\@gtempa=\@{}\the\@temptokena}}
\def\ps@pprintTitle{%
     \let\@oddhead\@empty\let\@evenhead\@empty%
     \let\@oddfoot\@empty\let\@evenfoot\@oddfoot}
\makeatother
\newtheorem{theo}{Theorem}[section]
\newtheorem{prop}[theo]{Proposition}

\newtheorem{lemm}[theo]{Lemma}
\newtheorem{exam}[theo]{Example}
\newtheorem{problem}[theo]{Problem}
\newdefinition{defi}[theo]{Definition}
\newdefinition{rema}[theo]{Remark}
\newdefinition{nota}[theo]{Notation}
\newproof{demo}{Proof}
\DeclareRobustCommand{\QED}{%
 \leavevmode\unskip\penalty9999\hbox{}\nobreak\hfill\quad\hbox{\qed}}
\newenvironment{smatrix}{\left[\begin{smallmatrix}}{\end{smallmatrix}\right]}
\def\Span<#1>{\langle#1\rangle}
\def\Tdef#1{\textbf{#1}}
\def\funct#1(#2)=#3:#4\to#5.{
  \begin{array}{r@{\,}c@{\,}c@{\,}l}#1:&#4&\to&#5\\&#2&\mapsto&#3\end{array}}
\MathOps{col=\vv,img,rk=rank,diag}
\MathSyms{bb}[v]{C,E,F,N,Z,R,v}
\MathSyms{cal}[c]{V,I,P,M,E,L,R}
\MathSyms{scr}[s]{S,V}
\def\Ev#1{\cE_{#1}}
\def\El#1{\cL_{#1}}
\def\Er#1{\cR_{#1}}
\def\PowR#1^#2{{\overline #1}^#2}

\let\imp\Rightarrow     \let\pmi\Leftarrow
\let\iff\Leftrightarrow 
\let\tens\otimes        \def\T{^\top\!}

\def\ABS{_{A,B;S}}

\begin{document}
\begin{frontmatter}
  \title{Linear spanning sets for matrix spaces}
  \author[uzh]{G.~Micheli\fnref{fn1}} \ead{giacomo.micheli at math.uzh.ch}
  \author[uzh]{J.~Rosenthal\fnref{fn1}} \ead{rosenthal at math.uzh.ch}
  \author[ua]{P.~Vettori\fnref{fn1,fn2}\corref{cor}} \ead{pvettori at
    ua.pt} \cortext[cor]{Corresponding author}
    \fntext[fn1]{%
     Authors supported in part by Swiss National Science Foundation grant
    SNF no. 149716.}
    \fntext[fn2]{%
This work was supported by Portuguese funds through the CIDMA (Center for Research and Development in Mathematics and Applications) and the Portuguese Foundation for Science and Technology (``FCT--Fundação para a Ciência e a Tecnologia''), within project PEst-OE/MAT/UI4106/2014.}
  \address[uzh]{University of Zurich, Winterthurstrasse 190, CH-8057 Zürich, Switzerland}
  \address[ua]{University of Aveiro, Campus de Santiago, 3810-193 Aveiro, Portugal}
  \begin{abstract}
    Necessary and sufficient conditions are given on matrices $A$, $B$
    and $S$, having entries in some field $\vF$ and suitable
    dimensions, such that the linear span of the terms $A^iSB^j$ over
    $\vF$ is equal to the whole matrix space.

    This result is then used to determine the cardinality of subsets of $\vF[A]S\vF[B]$ when $\vF$ is a finite field.
  \end{abstract}

\end{frontmatter}

\section{Introduction}

We start by stating a purely linear algebra problem:
\begin{problem}\label{prob}
  Let $m,n$ be integers and $\vF$ be any field. Let $A,S,B$ be
  matrices having entries in $\vF$ of dimensions $m\times m$, $m\times
  n$ and $n\times n$ respectively. Give necessary and sufficient
  conditions for the $\vF$-linear span of $\{A^i S B^j\}_{i,j\in\vN}$
  to be equal to the whole matrix space $\vF^{m\times n}$.
\end{problem}
A solution to this problem will be provided in Section \ref{basissec}.

Starting with Section~\ref{irredsec} we will assume that the base
field $\vF$ represents the finite field $\vF=\vF_q$ having cardinality
$q$.  Under these conditions and the conditions that $\gcd(m,n)=1$ and
the characteristic polynomials of the matrices $A$ and $B$ are
irreducible we are able to show in Section~\ref{irredsec} that $\{A^i
S B^j\}_{i,j\in\vN}$ spans the whole vector space $\vF^{m\times n}$ as
soon as $S\neq 0$.

In Section \ref{cardsec} we will prove that whenever the set $\{A^i S
B^j\}_{i,j\in\vN}$ spans the whole matrix ring as a vector space over
the finite field $\vF$, we are able to explicitly compute the
cardinality $\vF[A]S\vF[B]$.  A particular instance of this
computation (i.e. when $S$ is the identity matrix and $A$, $B$ have
irreducible characteristic polynomial) has already been approached via
inequalities in \cite{Chang13}.

\section{Notation and Preliminaries}

Let $\vF$ be a field and denote by $\Span<\sS>_\vF$ the linear span over $\vF$ of a set $\sS$ of elements in some $\vF$-vector space.
Entries, rows and columns of matrices are indexed by integers
starting from zero; $I_n$ and, respectively, $0_{m\times n}$ denote
the $n\times n$ identity matrix and the $m\times n$ zero matrix ---
indices may be omitted when no ambiguity arises.

Moreover, given $M\in\vF^{n\times n}$,
\begin{itemize}
\item the \Tdef{minimal polynomial} $\mu_M$ of $M$ is the monic
  generator of the ideal $\{p(s)\in\vF[s]:p(M)=0\}$;
\item the \Tdef{characteristic polynomial} of $M$ is
  $\chi_M(s)=\det(sI-M)$;
\item $\Ev M$ is the set of eigenvalues of $M$, i.e., the zeros of
  $\chi_M$ in some field extension of $\vF$;
\item $\El{M}^\lambda$ and $\Er{M}^\lambda$ are the left and,
  respectively, right eigenspaces of $M$ associated with
  $\lambda\in\Ev M$;
\item $\El{M}=\bigcup\limits_{\lambda\in\Ev
    M}\El{M}^\lambda\setminus\{0\}$ and
  $\Er{M}=\bigcup\limits_{\lambda\in\Ev
    M}\Er{M}^\lambda\setminus\{0\}$ are the sets of left and,
  respectively, right eigenvectors of $M$.
\item $M$ is \Tdef{cyclic} (or non-derogatory) if one of the following
  equivalent conditions holds true:
  \begin{itemize}
  \item $\mu_M=\chi_M$;
  \item $M$ is similar to a companion matrix;
  \item each eigenspace of $M$ has dimension 1, i.e., every
    eigenvector has geometric multiplicity 1.
  \end{itemize}
\end{itemize}

The definition of the Kronecker product and some of its properties are
given next. More details may be found in~\cite[Section
12.1]{LancTism85}.

\begin{defi}
  The \Tdef{Kronecker product} of matrices $M\in\vF^{m\times p}$ and
  $N\in\vF^{n\times q}$ is the block matrix
  \begin{equation*}
    M\tens N=[m_{i,j}N]_{0\le i<m,0\le j<p}\in\vF^{mn\times pq},
  \end{equation*}
  representing the tensor product of the linear maps corresponding to
  $M$ and $N$. Therefore, it satisfies the property
  \begin{equation}\label{eq.tpp}
    (M\tens N)(P\tens Q)=MP\tens NQ,
  \end{equation}
  whenever the matrix products on the right side can be computed.

  The \Tdef{(column) vectorization} of $M$ is the (column) vector
  $\col(M)\in\vF^{mp}$ formed by stacking the columns of $M$.
  Note that $\col:\vF^{m\times p}\to\vF^{mp}$ is an isomorphism of
  $\vF$-vector spaces, establishing a correspondence between
  entry $(i,j)$ of $M$ and entry $i+mj$ of $\col(M)$.

\end{defi}

Using this notation, given three matrices $M,X,N$ of suitable
dimensions,
\begin{equation}\label{eq.cmxn}
  \col(MXN)=(N\T\tens M)\col(X).
\end{equation}

\section{A basis for the vector space of $m\times n$ matrices}
\label{basissec}

Let matrices $A$, $B$, and $S$ as in Problem~\ref{prob} and define
\[
\sV\ABS=\Span<\{A^iSB^j\}_{i,j\ge0}>_\vF.
\]
In this and in the following section, conditions will be given that
ensure that the dimension of $\sV\ABS$ is maximal, i.e., equal to
$mn$.

\begin{theo}\label{th.mnge}
  Let $A\in\vF^{m\times m}$, $B\in\vF^{n\times n}$, and
  $S\in\vF^{m\times n}$ and consider the following conditions:
  \begin{align}
    \label{co.mngea}&
    \sV\ABS=\vF^{m\times n};\\
    \label{co.mngeb}&
    \text{$A$ and $B$ are cyclic};\\
    \label{co.mngec}& uSv\ne0,\;\forall u\in\El A,v\in\Er B.
  \end{align}
  Then, \eqref{co.mngea} $\iff$ $\big($\eqref{co.mngeb} and
  \eqref{co.mngec}$\big)$.
\end{theo}

\begin{rema}
  The previous theorem has also an impact in Cryptography since it
  gives necessary and sufficient conditions for the attack in
  \cite[Section 3]{Micheli13} to be performed in \emph{provable}
  polynomial time.
\end{rema}

Before proving the theorem, two lemmas will be stated.  The first one
provides a logical equivalence, which will be used within different
proofs.

\begin{lemm}\label{le.leq}
  Given three conditions $A$, $B$, and $C$, then $A\iff(B\text{ and
  }C)$ is equivalent to: $(A\imp B)$ and $\big(B\imp(A\iff C)\big)$.
\end{lemm}
\begin{demo}
  It is easy to check that both conditions are equivalent to the negation of $A$, when $B$ is false, and to $A\iff C$, when $B$ is true.
\end{demo}

The second lemma is well known (see~\cite{Hautus69,Shemesh84}) in the
case $\vF=\vC$. For completeness, a self-contained proof will be given here.
\begin{lemm}\label{le.pbh}
  Let $H\in\vF^{p\times p}$, $K\in\vF^{p\times q}$ and assume that
  $\Ev H\subseteq\vE$, extension field of $\vF$. Then, for any
  $d\ge\deg \mu_H$.
  \begin{equation*}
    \rk_\vF\begin{bmatrix}K&HK&\cdots&H^{d-1}K\end{bmatrix}=p\iff
    \rk_\vE\begin{bmatrix}\lambda I-H&K\end{bmatrix}=p,\;\forall \lambda\in\Ev H.
  \end{equation*}
\end{lemm}
\begin{demo}
  Observe that for any matrix $M$ with entries in $\vF$, $\rk_\vF
  M=\rk_\vE M$, since the rank depends only on the invertibility (in
  $\vF$) of square submatrices of $M$. So, this equivalent statement
  will be proved:
  \begin{equation*}
    \rk_\vE\begin{bmatrix}K&HK&\cdots&H^{d-1}K\end{bmatrix}<p\iff
    \exists\lambda\in\Ev H:\rk_\vE\begin{bmatrix}sI-H&K\end{bmatrix}<p.
  \end{equation*}
  ``$\imp$'': Be $u\in\vE^{1\times p}$ a nonzero vector such that
  $u\begin{bmatrix}K&HK&\cdots&H^{d-1}K\end{bmatrix}=0$ and be
  $a\in\vE[s]$ any generator of the principal ideal
  $\cI=\{f\in\vE[s]:uf(H)=0\}$. Since $\mu_H\in\cI$, $\deg
  a\le\deg\mu_H\le d$ and $a(\lambda)=0$ for some $\lambda\in\Ev H$.
  Write $a(s)=(\lambda-s)b(s)$, being $b(s)=\sum_{i=0}^{d-1}
  b_is^i\not\in\cI$. Hence, $v=ub(H)\ne0$. Moreover,
  \begin{equation*}
    vK=ub(H)K=
    \sum_{i=0}^{d-1} b_iuH^iK=
    \sum_{i=0}^{d-1} b_i0=0
  \end{equation*}
  and
  $0=ua(H)=u(\lambda I-H)b(H)=v(\lambda I-H)$. Thus,
  $v\begin{bmatrix}\lambda I-H&K\end{bmatrix}=0$.

  ``$\pmi$'':
  There exist $\lambda\in\Ev H$ and a nonzero $u\in\vE^{1\times p}$
  such that $u\begin{bmatrix}\lambda I-H&K\end{bmatrix}=0$, i.e.,
  $uH=\lambda u$ and $uK=0$. Hence,
  \begin{equation*}
    u\begin{bmatrix}K&HK&\cdots&H^{d-1}K\end{bmatrix}=
    u\begin{bmatrix}K&\lambda K&\cdots&\lambda^{d-1}K\end{bmatrix}=0.
  \end{equation*}
\end{demo}

\begin{demo}[of Theorem~\ref{th.mnge}]\def\refmnge#1{(\ref{co.mngeb}#1)}
  Consider the new conditions~\refmnge a: $A$ is cyclic and~\refmnge
  b: $B$ is cyclic, so that~\eqref{co.mngeb} is equivalent to~\refmnge
  a and~\refmnge b. Therefore, the equivalence
  $\eqref{co.mngea}\iff\big(\refmnge a\text{ and }\refmnge b\text{ and
  }\eqref{co.mngec}\big)$ will be proved.

  First of all, note that matrices $\{A^iSB^j\}$ generate
  $\vF^{m\times n}$ if and only if the corresponding vectors
  $\{\col(A^iSB^j)\}$ generate $\vF^{mn}$. Therefore, we get that
  \begin{equation}
    \eqref{co.mngea}\iff\label{eq.qed1}
    \Span<\{\col(A^iSB^j)\}_{i,j\ge0}>_\vF=
    \vF^{mn}.
  \end{equation}

  By~\eqref{eq.cmxn} and~\eqref{eq.tpp}, it follows that
  \begin{equation*}
    \col(A^iSB^j)=
    \col(A^iSB^jI_n)=
    (I_n\tens A^i)\col(SB^j)=
    (I_n\tens A)^i\col(SB^j).
  \end{equation*}
  Let $F=I_n\tens A\in\vF^{mn\times mn}$, which is a block diagonal
  matrix, and be $G$ the $mn\times n$ matrix whose columns are
  $\col(SB^j)$, $0\le j<n$.
  The (right) image of $G$, i.e., its column span, corresponds through $\col$ to the span of $SB^j$, $0\le j<n$. Analogously, for any $0\le i<m$, the image of $F^iG$ corresponds to the span of $A^iSB^j$, $0\le j<n$. Hence, by the Cayley-Hamilton Theorem,
  \begin{equation}\label{eq.qed2}
    \Span<\{\col(A^iSB^j)\}_{i,j\ge0}>_\vF
    \!=
    \img_\vF\begin{bmatrix}
      G\!&\!\!FG\!\!&\!\cdots\!&\!\!F^{m-1}G
    \end{bmatrix}.
  \end{equation}
  Observe that the degree of the minimal polynomial $\mu_F=\mu_{I\tens
    A}=\mu_A$ cannot be greater than $m$ and so, by~\eqref{eq.qed1},
  \eqref{eq.qed2} and Lemma~\ref{le.pbh}, we can state that
  \begin{align}\nonumber
    \eqref{co.mngea}&\iff\img_\vF\begin{bmatrix} G&FG&\cdots&F^{m-1}G
    \end{bmatrix}=\vF^{mn}\\
    &\iff \rk_\vE\begin{bmatrix} \lambda
      I-F&G\end{bmatrix}=mn,\;\forall\lambda\in\Ev A,\label{co.mnged}
  \end{align}
  being $\vE$ the extension field of $\vF$ containing the eigenvalues of $F$,
  i.e., of $A$.

  In order to determine the conditions that guarantee that the rank of
  the polynomial matrix $C(s)=\begin{bmatrix}sI-F&G\end{bmatrix}$ does
  not drop as $s\in\Ev A$, it is necessary to analyze the structure of $C(s)$ with greater detail.

  Denote by $G_i$, $0\le i<n$, the $m\times n$ blocks forming matrix
  $G$. Then
  \begin{align}
    C(s)=
    \begin{bmatrix}sI-F&G\end{bmatrix}=
    \begin{bmatrix}
      sI-A&&&&G_0\\
      &sI-A&&&G_1\\
      &&\ddots&&\vdots\\
      &&&sI-A&G_{n-1}
    \end{bmatrix}.
  \end{align}

  Now, let $\alpha$ be any eigenvalue of $A$ with geometric
  multiplicity $h$ and observe that the rank of the block-diagonal
  matrix $\alpha I-F$ (the first $mn$ columns of $C(\alpha)$) is equal
  to $n(m-h)=mn-nh$. Since matrix $G$ has $n$ columns, the rank of
  $C(\alpha)$ cannot exceed $mn-nh+n=mn-n(h-1)$. This shows that for
  condition~\eqref{co.mnged} to hold, it is necessary to have $h=1$,
  i.e., $A$ must be cyclic --- by equivalence~\eqref{co.mnged}, this
  shows that $\eqref{co.mngea}\imp\refmnge a$.

  On the other hand, by assuming that $A$ is cyclic, it follows that
  the rank of $C(\alpha)$ is $mn$ if and only if for every $w\ne0$
  such that $w(\alpha I-F)=0$, we have that $wC(\alpha)\ne0$. Since
  $\alpha I-F=I_n\tens(\alpha I-A)$, it turns out that $w(\alpha
  I-F)=0$ if and only if
  $w=\begin{bmatrix}u_0&u_1&\cdots&u_{n-1}\end{bmatrix}$, with
  $u_i\in\El A^\alpha$, $0\le i<n$. Therefore,
  \begin{align}\nonumber
    wC(\alpha)&=\begin{bmatrix}u_0&u_1&\cdots&u_{n-1}\end{bmatrix}
    \begin{bmatrix}
      \alpha I-A&&&&G_0\\
      &\alpha I-A&&&G_1\\
      &&\ddots&&\vdots\\
      &&&\alpha I-A&G_{n-1}
    \end{bmatrix}\\&=
    \begin{bmatrix}
      0&u_0G_0+u_1G_1+\cdots+u_{n-1}G_{n-1}
    \end{bmatrix}=
    \begin{bmatrix}0&g\end{bmatrix},\; g\in\vE^{1\times
      n}.\label{eq.wcm}
  \end{align}

  Since the eigenspace $\El A^\alpha$ has dimension $1$, is it
  generated by one (eigen)vector, say $u\ne0$, whence $u_i=\gamma_i u$,
  $\gamma_i\in\vE$ for $0\le i<n$, not all zero.  This means that
  \[g=\gamma_0uG_0+\gamma_1uG_1+\cdots+\gamma_{n-1}uG_{n-1}\] is not
  zero if and only if vectors $\{uG_i\}_{0\le i<n}$ are linearly
  independent. Hence, by equivalence~\eqref{co.mnged},
  condition~\refmnge a implies that $\eqref{co.mngea}$ is equivalent
  to the linear independence of $\{uG_i\}_{0\le i<n}$, for every
  $u\in\El A$. We already proved that $\eqref{co.mngea}\imp\refmnge a$
  and so, by Lemma~\ref{le.leq}, it follows that
  \begin{equation}\label{eq.qed4}
    \eqref{co.mngea}\iff\big(\refmnge a\text{ and }\forall u\in\El A,
    \{uG_i\}_{0\le i<n}\text{ are $\vE$-linearly independent}\big).
  \end{equation}

  Consider now any $u\in\vE^{1\times m}$ and define the matrix
  \begin{equation*}
    D=(I_n\tens u)G=
    \begin{bmatrix}
      uG_0\\uG_1\\\vdots\\uG_{n-1}
    \end{bmatrix}\in\vE^{n\times n}.
  \end{equation*}
  Moreover, for every $0\le i<n$ and $0\le j<n$, let $(SB^j)_i$ be the $i$-th column of $SB^j$.

  By definition, the $j$-th column of $G$ is $\col(SB^j)$, which
  contains, stacked, vectors $(SB^j)_i$. Therefore, in particular,
  the $j$-th column of $G_i$, is $(SB^j)_i$. Consequently, the $j$-th
  component of $uG_i$, which is the entry at $(i,j)$ of $D$, is
  $u(SB^j)_i$.  At the same time, this value is the $i$-th component
  (column) of $uSB^j$, i.e, the entry at $(j,i)$ of the matrix whose
  rows are $uSB^j$. In other words,
  \begin{equation*}
    D\T=\begin{bmatrix}
      uSB^0\\uSB^1\\\vdots\\uSB^{n-1}
    \end{bmatrix}.
  \end{equation*}
  Since $D$ is square, its rows are linearly independent if and only
  if its columns share the same property.
  Applying again Lemma~\ref{le.pbh} with $H=B\T$ and $K=(uS)\T$, we
  get that
  \begin{align}
    &\{uG_i:0\le i<n\} \text{ are $\vE$-linearly independent }
    \label{eq.qed6}&\iff\\
    &\{uSB^j:0\le j<n\} \text{ are $\vE$-linearly independent }
    \nonumber&\iff\\
    &\rk_\vE\begin{bmatrix}(uS)\T&B\T(uS)\T&\cdots&(B\T)^{n-1}(uS)\T\end{bmatrix}=n
    \nonumber&\iff\\
    &\rk_\vE\begin{bmatrix}\lambda I-B\T&(uS)\T\end{bmatrix}=
    \rk_\vE\begin{bmatrix}\lambda I-B\\uS\end{bmatrix}=n,\;
    \forall\lambda\in\Ev B.\nonumber
  \end{align}
  As before, consider
  $E(s)=\left[\begin{smallmatrix}
      sI-B\\uS\end{smallmatrix}\right]\in\vE^{(n+1)\times n}[s]$
  and any $\beta\in\Ev B$.  Since $\beta I-B$ has rank $n-k$, where
  $k$ is the geometric multiplicity of $\beta$, $\rk_\vE E(\beta)\le
  n-k+1$. We conclude that, when~\eqref{eq.qed6} holds, then $k=1$,
  i.e., $\eqref{eq.qed6}\imp\refmnge b$, i.e., $B$ is cyclic.

  By assuming that $B$ is cyclic, the rank of $E(\beta)$ is
  effectively $n$ if $E(\beta)v\ne0$ for any $v\in\Er B^\beta$.
  Since $(\beta I-B)v=0$, condition $E(\beta)v\ne0$ reduces to
  $uSv\ne0$: \begin{equation*} \text{if $B$ is cyclic, i.e., \refmnge
      b holds, }\eqref{eq.qed6}\iff uSv\ne0,\;\forall v\in\Er{B}.
  \end{equation*}
  Thus, by Lemma~\ref{le.leq}, $\eqref{eq.qed6}\iff\big(\refmnge
  b\text{ and }uSv\ne0,\;\forall v\in\Er{B}\big)$.  This, together
  with~\eqref{eq.qed4}, concludes the proof.
\end{demo}

\begin{exam}
  Consider the following matrices, with $m,n\ge2$:
  \begin{align*}
    A&\!=\!\arraycolsep3pt
    \begin{bmatrix}0&0\\I_{m-1}&0\end{bmatrix}\!\in\!\vF^{m\times m},&
    B&\!=\!\arraycolsep3pt
    \begin{bmatrix}0&I_{n-1}\\0&0\end{bmatrix}\!\in\!\vF^{n\times n},&
    S&\!=\!\arraycolsep3pt
    \begin{bmatrix}1&0\\0&0_{(m-1)\times(n-1)}\end{bmatrix}\!\in\!\vF^{m\times
      n}.
  \end{align*}
  Both $A$ and $B$ are already in (left and right, respectively)
  Jordan canonical form. Therefore, their only eigenvalue is
  $\lambda=0$, they are nilpotent and cyclic with minimal polynomials
  $\mu_A(s)=s^m$ and $\mu_B(s)=s^n$, and their eigenspaces are
  generated by $u=\begin{bmatrix}1&0&\cdots&0\end{bmatrix}$ (left
  eigenvector of $A$) and
  $v=\begin{bmatrix}1&0&\cdots&0\end{bmatrix}\T$ (right eigenvector of
  $B$).

  Even though $S$ has rank 1, $uSv=1\ne0$, whence conditions~\eqref{co.mngeb} and~\eqref{co.mngec} of
  Theorem~\ref{th.mnge} are satisfied. Therefore, $\vF$--linear
  combinations of matrices $E_{i,j}=A^iSB^j$, with $0\le i<m$ and
  $0\le j<n$, generate $\vF^{m\times n}$ for any field $\vF$.

  Indeed, it is straightforward to check that each $E_{i,j}$ is one of
  the $mn$ elements of the canonical basis of $\vF^{m\times n}$,
  having its unique nonzero entry, equal to $1$, at position
  $(i,j)$.
  In other words, $\col(E_{i,j})$ is the $i+mj$-th vector of the
  canonical basis of $\vF^{mn}$.
\end{exam}

To the authors' knowledge, equality~\eqref{co.mngea} and the kind of
equivalent conditions that were presented in Theorem~\ref{th.mnge}
have not been considered in the literature before (not even when
$m=n$: see, for instance, the survey~\cite{Laffey86} containing a
small section about spanning sets of matrix algebras).

A comparison with previous results can be made only in the case
$m=n=2$, verifying that $\vF^{2\times 2}$ is spanned by linear
combinations of $A^iB^j$, $i,j=0,1$, if and only if it can be
generated by $A$ and $B$ as a matrix algebra. (The well-known
criterium for the latter problem, presented in the following
proposition, can be found, for example, in~\cite{AslaSlet09}, where it
is thoroughly investigated.)

\begin{prop}\label{pr.tbt}
  Let $A,B\in\vF^{2\times2}$ and $S=I$. Then, the commutator
  $[A,B]=AB-BA$ is invertible if and only if
  conditions~\eqref{co.mngeb} and~\eqref{co.mngec} hold.
\end{prop}
\begin{demo}
  Notice that adding a scalar matrix $cI$, $c\in\vF$, to $A$ or $B$
  does not change both the spanned space and the generated algebra, nor the commutator $[A,B]$.
  Therefore, we shall assume that $A$ and $B$ have zero trace.

  First, observe that $A$ is not cyclic if and only if its canonical
  Jordan form is a scalar matrix if and only if $A$ itself is a scalar
  matrix, i.e., zero. Therefore, if either $A$ or $B$ is not cyclic,
  $[A,B]=0$. This proves that
  \begin{equation} [A,B]\text{ is invertible }\imp\text{
      \eqref{co.mngeb}}.
    \label{eq.ab1}
  \end{equation}

  Assume now~\eqref{co.mngeb}, both $A$ and $B$ are cyclic, and
  suppose, without loss of generality, that $A$ is in Jordan form.
  This means that
  \begin{equation}\label{eq.ab2}
    A=\begin{bmatrix}a&b\\0&-a\end{bmatrix},\quad
    B=\begin{bmatrix}\alpha&\beta\\\gamma&-\alpha\end{bmatrix},\quad
    \text{and}\quad
    [A,B]=\begin{bmatrix}b\gamma&2a\beta-2\alpha b\\-2a\gamma&-b\gamma\end{bmatrix}.
  \end{equation}
  In order to be cyclic, i.e., not zero, matrix $B$ must satisfy $\alpha\ne0$, $\beta\ne0$ or $\gamma\ne0$.
  For matrix $A$, the two following cases are possible.
  \begin{enumerate}
  \item$a=0$ and $b=1$:
    $\El A\cup\{0\}=\El A^0$ is generated by
    $u=\begin{bmatrix}0&1\end{bmatrix}$.  If $\gamma=0$ then
    $\alpha\in\Ev B$ and $v=\begin{smatrix}1\\0\end{smatrix}\in\Er
    B^\alpha$, satisfying $uv=0$. On the other hand, if $uv=0$, with
    $v\in\Er B^\lambda$ for some $\lambda\in\Ev B$, then
    $v=\begin{smatrix}x\\0\end{smatrix}$, $x\ne0$. By definition,
    \[
    Bv=\lambda v\iff
    \begin{bmatrix}\alpha x\\\gamma x\end{bmatrix}=
    \begin{bmatrix}\lambda x\\0\end{bmatrix}\iff
    \gamma=0.
    \]
    By~\eqref{eq.ab2}, it easy to check that $[A,B]$ is singular if
    and only if $\gamma=0$, thus proving that $[A,B]$ invertible
    $\iff$ \eqref{co.mngec}.

  \item$a\ne0$ and $b=0$: both $u=\begin{bmatrix}1&0\end{bmatrix}$ and
    $u=\begin{bmatrix}0&1\end{bmatrix}$ belong to $\El A$.
    If $\beta\gamma=0$ then, similarly to the previous case,
    $\begin{smatrix}1\\0\end{smatrix}$ or
    $\begin{smatrix}0\\1\end{smatrix}$ belong to $\Er B$, being
    possible to satisfy $uv=0$ with a nonzero $v\in\Er B$.
    Vice versa, if $uv=0$ for some
    $v=\begin{smatrix}x\\y\end{smatrix}\in\Er B$, then $xy=0$. It
    turns out that $Bv=\lambda v$ implies that $x=0\imp\beta=0$ and
    $y=0\imp\gamma=0$, therefore $\beta\gamma=0$.

    Concluding, by~\eqref{eq.ab2}, $[A,B]$ invertible $\iff$
    $\beta\gamma\ne0$ $\iff$ \eqref{co.mngec}.
  \end{enumerate}
  We showed that, in both cases, when~\eqref{co.mngeb} holds, then
  $[A,B]$ is invertible $\iff$ \eqref{co.mngec}.  The statement
  follows by~\eqref{eq.ab1} and Lemma~\ref{le.leq}.
\end{demo}

When conditions~\eqref{co.mngeb} and~\eqref{co.mngec} of
Theorem~\ref{th.mnge} are not satisfied, matrices $A^iSB^j$, with
$0\le i<m$ and $0\le j<n$, are linearly dependent. However, something more can
be said about the dimension of the space they generate.

The general case demands an extremely complicated notation: only the case of cyclic and diagonalizable matrices $A$ and $B$ will be considered in this paper.

\begin{theo}\label{th.gdsm}
  Let $S\in\vF^{m\times n}$ and suppose that $A\in\vF^{m\times m}$ and
  $B\in\vF^{n\times n}$ are cyclic and diagonalizable. In particular,
  be $U\in\vE^{m\times m}$ and $V\in\vE^{n\times n}$ two invertible
  matrices, in some extension field $\vE$ of $\vF$, such that
  $UAU^{-1}$ and $V^{-1}BV$ are diagonal.

  Then, the dimension of $\sV\ABS$, is equal to the number of nonzero
  entries of $USV$.
\end{theo}

Before proving Theorem~\ref{th.gdsm}, we introduce the necessary
notation and state a fundamental lemma.

Given $A\in\vF^{m\times m}$, $B\in\vF^{n\times n}$, and
$S\in\vF^{m\times n}$, let $r_{i,j}=\col(A^iSB^j)$ and define
\begin{equation}\label{eq.drabs}
  R\ABS=
  \begin{bmatrix}
    r_{0,0}\!&\!  r_{1,0}\!&\!  \cdots\!&\!  r_{m-1,0}\!&\!
    r_{0,1}\!&\!  r_{1,1}\!&\!  \cdots\!&\!  r_{m-1,n-1}
  \end{bmatrix}\in\vF^{mn\times mn}.
\end{equation}

Then, given $v\in\vF^n$, $\diag(v)\in\vF^{n\times n}$ is the diagonal
matrix defined by the components of $v$. Moreover, let
$\diag(M)=\diag\big(\col(M)\big)$ for any matrix $M$.

Finally, let $\PowR
x^n=\begin{bmatrix}1&x&\cdots&x^{n-1}\end{bmatrix}$ and be
$\cV^n_{x_1,\ldots,x_k}$ the matrix whose rows are $\PowR x^n_1$,
\ldots, $\PowR x^n_k$.

\begin{lemm}\label{le.rabs}
  Let $A\in\vF^{m\times m}$, $B\in\vF^{n\times n}$, and
  $S\in\vF^{m\times n}$. Suppose that $u_h\in\El A^{\alpha_h}$, $0\le
  h<s$, and $v_k\in\El B^{\beta_k}$, $0\le k<t$, are the rows and,
  respectively, columns of matrices $U\in\vE^{s\times m}$ and
  $V\in\vE^{n\times t}$ in a suitable extension field $\vE$ of $\vF$. Then,
  \begin{equation}\label{eq.rabs}
    (V\T\tens U)R\ABS=\diag(USV)
    (\cV^n_{\beta_1,\ldots,\beta_t}\tens\cV^m_{\alpha_1,\ldots,\alpha_s}).
  \end{equation}
\end{lemm}

\begin{demo}
  Observe that, for any row $u_h$ of $U$ and column $v_k$ of $V$,
  there exist $\alpha_h\in\Ev A$ and $\beta_k\in\Ev B$ such that
  $u_h\in\El A^{\alpha_h}$ and $v_k\in\Er B^{\beta_k}$.  Thus,
  \begin{equation*}
    (v_k\T\tens u_h)\col(A^iSB^j)=
    u_hA^iSB^jv_k=
    u_hSv_k\,\alpha_h^i\beta_k^j
  \end{equation*}
  and, from~\eqref{eq.drabs}, it follows that
  \begin{equation*}
    (v_k\T\tens u_h)R\ABS=
    u_hSv_k\big(\PowR\beta^n_k\tens\PowR\alpha^m_h\big).
  \end{equation*}
  \emph{Stacking up} all these equalities, we get
  equation~\eqref{eq.rabs}.
\end{demo}

\begin{rema}
  Using Lemma~\ref{le.rabs}, implication
  $\eqref{co.mngea}\imp\eqref{co.mngec}$ of Theorem~\ref{th.mnge} can
  be proved in a much simpler way.

  Indeed, suppose that the nonzero left-eigenvector $u\in\El A^\alpha$
  and right-eigenvector $v\in\Er B^\beta$ satisfy $uSv=0$. Then,
  taking $U=u$ and $V=v$ in formula~\eqref{eq.rabs}, we get
  \[
  (v\T\tens u)R\ABS=(uSv)(\PowR \beta^n\tens\PowR\alpha^m)=0,
  \]
  showing that $R\ABS$ does not have full rank. Therefore, its columns
  $\col(A^iSB^j)$ are linearly dependent and the set of matrices
  $A^iSB^j$ cannot generate $\vF^{m\times n}$.
\end{rema}

\begin{demo}[of Theorem~\ref{th.gdsm}]
  Let $\alpha_h$, $0\le h<m$ and $\beta_k$, $0\le k<n$, be the left
  eigenvalues of $A$ associated with the rows of $U$ and,
  respectively, the right eigenvalues of $B$ associated with the
  columns of $V.$

  Since $A$ and $B$ are cyclic and diagonalizable, they have no
  repeated eigenvalues, whence $\cV^m_{\alpha_1,\ldots,\alpha_s}$ and
  $\cV^n_{\beta_1,\ldots,\beta_n}$ are invertible Vandermonde
  matrices.

  By Lemma~\ref{le.rabs}, we have that
  \[
  (V\T\tens U)R\ABS=\diag(USV)
  (\cV^n_{\beta_1,\ldots,\beta_n}\tens\cV^m_{\alpha_1,\ldots,\alpha_s}),
  \]
  where both Kronecker products are invertible. So, $\rk R\ABS=\rk
  \diag(USV)$, which is equal to the number of nonzero entries of
  $USV$.

  Since by definition~\eqref{eq.rabs}, the (column) rank of $R\ABS$ is
  equal to the dimension of the space spanned by $A^iSB^j$, the proof
  is concluded.
\end{demo}

\section{The irreducible case}
\label{irredsec}

For the remainder of the paper we will asssume that $\vF=\vF_q$
represents the finite field of order $q$.

The main result of this section will provide a necessary and
sufficient condition for matrices $A$, $B$ having irreducible
characteristic polynomial which guarantees that
condition~\eqref{co.mngea} of Theorem~\ref{th.mnge} holds true:
\begin{theo}\label{th.abidm}
  Let $\vF$ be a finite field, $A\in\vF^{m\times m}$, $S\in
  \vF^{m\times n}$ and $B\in\vF^{n\times n}$.  Suppose that $A$ and
  $B$ have irreducible characteristic polynomials. Then,
  \[
  \sV\ABS=\vF^{m\times n},\forall S\ne0\text{ if and only if }\gcd(m,n)=1.
  \]
\end{theo}

\begin{demo}
  Define the $\vF$-linear map
  \begin{equation}\label{def.psi}
    \funct\psi(Z=[z_{i,j}])=\sum_{\substack{0\le i<m\\0\le j<n}}
    z_{i,j}A^iSB^j:\vF^{m\times n}\to\vF^{m\times n}.
  \end{equation}
  and note that $\sV\ABS$ is the image of $\psi$. Therefore, we need
  to prove that $\ker\psi=\{0\},\forall S\ne0\iff\gcd(m,n)=1$. By~\eqref{eq.cmxn} we
  obtain that
  \begin{align*}
    \col\big(\psi(Z)\big)&=
    \col\left(\sum_{\substack{0\le i<m\\0\le j<n}}z_{i,j}A^iSB^j\right)=
    \sum_{\substack{0\le i<m\\0\le j<n}}z_{i,j}(B^j)\T\tens A^i\col(S).
  \end{align*}
  Hence, by injectivity of $\col$, it follows that $\psi$ is injective
  (for any choice of $S\ne0$) if and only if the kernel of matrix
  $M=\sum_{0\le i<m,0\le j<n}z_{i,j}(B^j)\T\tens A^i$ is trivial, i.e., $M$ has
  no zero eigenvalues whenever $Z\ne0$.

  Observe first that, by the assumptions on $A$ and $B$, the matrix
  rings $\vF[A]$ and $\vF[B]$ are fields.  Moreover, all eigenvalue
  $\alpha\in\Ev A$ and $\beta\in\Ev B$ have $\vF$-linearly independent
  powers up to degree $m-1$ and, respectively, $n-1$, being
  $\vF(\alpha)\cong\vF[A]$ and $\vF(\beta)\cong\vF[B]$, which are
  Galois extensions of $\vF$ of degree $m$ and, respectively, $n$.

  By a classical result on Kronecker products (see,
  e.g.,~\cite[Theorem~1, p.~411]{LancTism85} for $\vF=\vR$, whose
  generalization to finite fields is straightorward) the set of
  eigenvalues of $M$ is
  \begin{equation}\label{eq.evtp}
    \Ev M=\left\{\sum_{\substack{0\le i<m\\0\le j<n}}z_{i,j}\alpha^i\beta^j:
      \alpha\in\Ev A,\beta\in\Ev B \right\},
  \end{equation}
  where all eigenvalues are considered as elements in some common
  field extension.

  So, $\ker\psi=\{0\}$ if and only if each sum in~\eqref{eq.evtp} is
  nonzero. In other words, for any two $\alpha\in\Ev A$ and
  $\beta\in\Ev B$, the products $\{\alpha^i\beta^j\}_{i<m,j<n}$ are
  $\vF$-linearly independent. By~\cite[Proposition~5.1 and
  Theorem~5.5]{Cohn91}, this condition is equivalent to
  \[\vF(\alpha)\cap\vF(\beta)=\vF.\]
  Since the intersection of $\vF(\alpha)$ and $\vF(\beta)$ is the
  field extension of $\vF$ of degree $\gcd(m,n)$
  (see~\cite[Theorem~2.6]{LidlNied97}), the proof is concluded.
\end{demo}

\section{The cardinality of subsets of $\vF[A]S\vF[B]$}
\label{cardsec}

In this section we will explicitly compute the cardinality of the set
$\vF[A]S\vF[B]$ whose relevance in Cryptography is discussed in
\cite{Chang13,MazeMoniRose07}.  Define the space of
polynomials
\[
\cP^k[s]=\{p(s)\in\vF[s]:\deg p<k\},\;k=0,1,\ldots
\]
being, for instance, $\cP^0=\{0\}$ and $\cP^1=\vF$.

Note that, given a square matrix $M$ with $d=\deg \mu_M$,
\[
\cP^0[M]\subset\cP^1[M]\subset\cdots\subset\cP^{d-1}[M]\subset\cP^d[M]=\cP^k[M],\;\forall
k\ge d.
\]

The main objective of this section consists in calculating the cardinality of
the set
\[
\cM\ABS^{h,k}=\cP^h[A]S\cP^k[B]\subseteq\vF^{m\times n}.
\]
\begin{theo}\label{th.cr}
  Let $A\in\vF^{m\times m}$, $B\in\vF^{n\times n}$, and
  $S\in\vF^{m\times n}$ such that $\sV\ABS=\vF^{m\times n}$. Then, for
  any $0\le h\le m$ and $0\le k\le n$,
  \[
  \left|\cM\ABS^{h,k}\right|=\frac{(q^h-1)(q^k-1)}{q-1}+1.
  \]
\end{theo}

In order to demonstrate this statement, some specific notation and one
preparatory lemma are needed.

First, for every $h\le m$, let
\[
\vF^{h;m}=\{x\in\vF^m:x_i=0,\forall i=h,\ldots,m-1\},
\]
being therefore $\vF^h\cong\vF^{h;m}\subseteq\vF^m$.  Define, for
every $h\le m$ and $k\le n$, the bilinear map
\begin{align}\label{def.phi}
  \funct\phi^{h,k}((x,y))=xy\T:\vF^{h;m}\times\vF^{k;n}\to\vF^{m\times
    n}.
\end{align}
and, for the sake of simplicity, denote its image by
\begin{align}\label{def.phk}
  \Phi^{h,k}=\phi^{h,k}(\vF^{h;m}\times\vF^{k;n}).
\end{align}

\begin{lemm}\label{lm.cr}
  Let $A$, $B$, and $S$ as in Theorem~\ref{th.cr}. Then
  $\left|\cM\ABS^{h,k}\right|=|\Phi^{h,k}|$.
\end{lemm}
\begin{demo}
  Consider the map $\psi$ defined in~\eqref{def.psi}.  We claim that
  $\psi(\Phi^{h,k})=\cM\ABS^{h,k}$. Actually, for every
  $M\in\cM\ABS^{h,k}$, there exist
  $(x,y)\in\vF^{h;m}\times\vF^{k;n}\subseteq\vF^m\times\vF^n$ such
  that
  \[
  M= \Bigg(\sum_{0\le i<h}x_iA^i\Bigg)S \Bigg(\sum_{0\le
    j<k}y_jB^j\Bigg)= \sum_{\substack{0\le i<m\\0\le
      j<n}}x_iy_jA^iSB^j=\psi(xy\T)\in\psi(\Phi^{h,k}).
  \]
  Therefore, $\left|\cM\ABS^{h,k}\right|\le|\Phi^{h,k}|$.  Moreover,
  when $\sV\ABS=\vF^{m\times n}$, $\psi$ is injective and so
  $\Phi^{h,k}\leftrightarrow\cM\ABS^{h,k}$.
\end{demo}

Observe that this lemma shows that the cardinality of $\cM\ABS^{h,k}$
is independent of the choice of $A$, $B$, and $S$ when
condition~\eqref{co.mngea} is met.

The problem is now reduced to the computation of the cardinality of $\Phi^{h,k}$,
defined in~\eqref{def.phk}.

\begin{demo}[of Theorem~\ref{th.cr}]
  Consider again the map $\phi^{h,k}$, defined in~\eqref{def.phi}, and
  observe that
  \[
  \vF^{h;m}\times\vF^{k;n}=(\phi^{h,k})^{-1}(\Phi^{h,k})=\bigcup_{Z\in\Phi^{h,k}}{(\phi^{h,k})^{-1}(Z)}.
  \]
  Consequently, since the inverse images are disjoint,
  \[
  q^hq^k= |\vF^{h;m}\times\vF^{k;n}|=
  \Bigg|\bigcup_{Z\in\Phi^{h,k}}(\phi^{h,k})^{-1}(Z)\Bigg|=
  \sum_{Z\in\Phi^{h,k}}|(\phi^{h,k})^{-1}(Z)|.
  \]
  To compute the value of the summation, we have to consider two
  situations.
  \begin{itemize}
  \item When $Z=0$, $\phi(x,y)=xy\T=0$ if and only if all the products
    of each component of $x$ and each component of $y$ are zero if and
    only if $x=0$ and $y=0$ ($1$ case), $x=0$ and $y\ne0$ ($q^k-1$
    cases), or $x\ne0$ and $y=0$ ($q^h-1$ cases). Therefore,
    $|\phi^{-1}(0)|= q^h+q^k-1$.
  \item If $Z\ne0$, observe that, by the bilinearity of $\phi^{h,k}$,
    $\phi^{h,k}(x,y)=\phi^{h,k}(\alpha x,\alpha^{-1}y)$ for every
    $\alpha\in\vF\setminus\{0\}$.

    On the other hand, if $\phi^{h,k}(x,y)=\phi^{h,k}(\tilde x,\tilde
    y)$ then $\tilde x=\alpha x$ and $\tilde y=\alpha^{-1}y$ for some
    $\alpha\ne0$. Indeed, considering only the indexes $i$ and $j$
    such that $x_iy_j=\tilde x_i\tilde y_j\ne0$, we get that
    \[\frac{x_i}{\tilde x_i}=\frac{\tilde y_j}{y_j}.\]
    By the independency of the indices, it follows that
    $\alpha=\frac{x_i}{\tilde x_i}=\frac{\tilde y_j}{y_j}$ for every
    $i,j$.  So, we conclude that $|(\phi^{h,k})^{-1}(Z)|=|\vF\setminus\{0\}|=q-1$.
  \end{itemize}
  Putting all together,
  \begin{align*}
    q^hq^k&= |(\phi^{h,k})^{-1}(0)|+
    \!\!\!\!\!\sum_{Z\in\Phi^{h,k}\setminus \{0\}}\!
    \left|(\phi^{h,k})^{-1}(Z)\right|\\&= q^h+q^k-1+
    \!\!\!\!\!\sum_{Z\in\Phi^{h,k}\setminus\{0\}}\!\!\!\!\!  (q-1)
    =q^h+q^k-1+\big(|\Phi^{h,k}|-1\big)(q-1),
  \end{align*}
  whence
  \[
  |\Phi^{h,k}|= \frac{q^hq^k-q^h-q^k+1}{q-1}+1=
  \frac{(q^h-1)(q^k-1)}{q-1}+1.
  \]
  Finally, the claim follows by Lemma~\ref{lm.cr}.
\end{demo}


\end{document}